\documentclass{colt2014} 
\usepackage[english]{babel}

\newcommand{\cA}{\mathcal{A}}
\newcommand{\cB}{\mathcal{B}}
\newcommand{\cC}{\mathcal{C}}

\newcommand{\cF}{\mathcal{F}}

\newcommand{\cM}{\mathcal{M}}
\newcommand{\cN}{\mathcal{N}}

\newcommand{\bE}{\mathbb{E}}
\newcommand{\bP}{\mathbb{P}}
\newcommand{\N}{\mathbb{N}}
\newcommand{\R}{\mathbb{R}}

\usepackage{dsfont}

\def \ind {\mathds{1}}
\newcommand{\norm}[2]{\mathcal{N}\left(#1,#2\right)}

\newcommand{\K}{\mathrm{KL}}
\newcommand{\Kb}{\mathrm{K}}

\newcommand{\argmax}[1]{\underset{#1}{\text{argmax}}}

\usepackage{algorithm,algorithmic}

\def \qed {\begin{flushright} $\qquad \Box $
  \end{flushright}}
  
\usepackage{blindtext}

\title[On the Complexity of A/B Testing]{On the Complexity of A/B Testing} 
\usepackage{times}

 \coltauthor{\Name{Emilie Kaufmann} \Email{kaufmann@telecom-paristech.fr}\\
 \addr LTCI, T\'el\'ecom ParisTech \& CNRS
 \AND
 \Name{Olivier Capp\'e} \Email{cappe@telecom-paristech.fr}\\
 \addr LTCI, T\'el\'ecom ParisTech \& CNRS
 \AND
 \Name{Aur\'elien Garivier} \Email{aurelien.garivier@math.univ-toulouse.fr}\\
 \addr Institut de Math\'ematiques de Toulouse, Universit\'e Paul Sabatier
 }

\begin{document}

\maketitle

\begin{abstract}
  A/B testing refers to the task of determining the best option among two
  alternatives that yield random outcomes. We provide
  distribution-dependent lower bounds for the performance of A/B
  testing that improve over the results currently available both in the
  fixed-confidence (or $\delta$-PAC) and fixed-budget settings. When the distribution of the outcomes  
  are Gaussian, we prove that the complexity of the fixed-confidence and
  fixed-budget settings are equivalent, and that uniform sampling of both
  alternatives is optimal only in the case of equal variances. In the common
  variance case, we also provide a stopping rule that terminates faster than
  existing fixed-confidence algorithms. In the case of Bernoulli distributions,
  we show that the complexity of fixed-budget setting is smaller than that
  of fixed-confidence setting and that uniform sampling of both
  alternatives---though not optimal---is advisable in practice when combined
  with an appropriate stopping criterion.
\end{abstract}


\section{Introduction}

A/B Testing is a popular procedure used, for instance, for website
optimization: two versions of a webpage, say A and B, are empirically
compared by being presented to users. Each user only sees one of the two
versions, and the goal is to determine which version is preferable. We assume
that the users provide a real-valued index of the quality of
the pages, which is modeled by probability distributions $\nu_A$ and
$\nu_B$, with respective means $\mu_A$ and $\mu_B$. For example, a standard
objective is to determine which webpage has the highest conversion
 rate (probability that a user actually becomes a customer) 
 by receiving binary feedback from the users.

Methods for A/B Testing are often viewed as statistical tests of the hypothesis
$H_0 : (\mu_A \leq \mu_B)$ against $H_1 : (\mu_A>\mu_B)$. One may consider either
\textit{classical tests}, based on a number of samples $n_A$ and $n_B$ from
each distribution fixed before the experiment, or
\textit{sequential tests}, based on paired samples ($X_s,Y_s$) of $\nu_A,\nu_B$
and in which a randomized stopping rule determines when the experiment is to be
terminated. In both of these test settings, the sampling schedule is determined in
advance, which is a possible source of sub-optimality as A/B Testing algorithms
could take advantage of past observations to provide a smarter choice of the page
to be displayed to the next user. In the sequel, we investigate whether A/B Testing 
could benefit from an adaptive sampling schedule. Ignoring the possible long-term effects on users of presenting one or the other option, we consider it as a particular instance of \textit{best arm identification in a two-armed bandit model}.

A two-armed bandit model consists of two unknown probability distributions on
$\R$, $\nu_1$ and $\nu_2$, sometimes referred to as arms or options (webpages in
our motivating example). Arm $a$ has expectation $\mu_a$.  At each time
$t=1,2,\dots$, an agent chooses an option $A_t\in\{1,2\}$ and receives an
independent draw $Z_t$ of the corresponding arm $\nu_{A_t}$.  We denote by
$\bP_\nu$ (resp. $\bE_\nu$) the probability law (resp. expectation) of the
corresponding process $(Z_t)$. We assume that the bandit model
$\nu=(\nu_1,\nu_2)$ belongs to a class $\cM$ such that for all $\nu\in\cM$,
$\mu_1\neq\mu_2$. 
In order to identify the best arm, that is the arm $a^*$ with highest
expectation, the agent must use a strategy defining which arms to sample from,
when to stop sampling, and which arm $\hat{a}$ to choose.  The \textit{sampling
  rule} determines how, at time $t$, the arm $A_t$ is chosen based on the past
observations; in other words, $A_t$ is $\cF_{t-1}$--measurable, with
$\mathcal{F}_t=\sigma(A_1,Z_1,\dots,A_t,Z_t)$. The \textit{stopping rule}
$\tau$ is a stopping time with respect to $(\mathcal{F}_{t})_{t\in\N}$
satisfying $\bP_\nu(\tau<+\infty)=1$.  The \textit{recommendation rule} is a
$\mathcal{F}_\tau$-measurable random arm $\hat{a} \in \{1,2\}$.  This triple
$((A_t),\tau,\hat{a})$ entirely determines the strategy, which we denote in the
sequel by $\cA$. As discussed before, statistical tests correspond to strategies that sample the arms in a round-robin fashion, which we will refer to as \textit{uniform sampling}.

In the bandit literature, two different settings have been considered.  In the
\emph{fixed-confidence setting}, a risk parameter $\delta$ is fixed.  A
strategy $\cA$ is called \emph{$\delta$-PAC} if, for every choice of
$\nu\in\mathcal{M}$, $\bP_\nu(\hat{a} = a^*)\geq 1-\delta$.  The goal is, among
the \emph{$\delta$-PAC} strategies, to minimize the expected number of draws
$\bE_\nu[\tau]$. In the \emph{fixed-budget setting}, the number of draws $\tau$
is fixed in advance ($\tau = t$ almost surely) and the goal is to choose the
recommendation rule so as to minimize $p_t(\nu) := \bP_\nu(\hat{a} \neq
a^*)$. In the fixed-budget setting, a strategy $\cA$ is called
\emph{consistent} if, for every choice of $\nu\in\mathcal{M}$, $p_t(\nu)
\underset{t\rightarrow \infty}{\rightarrow} 0.$

In order to unify and compare these approaches, we define the \emph{complexity}
$\kappa_C(\nu)$ (resp. $\kappa_B(\nu)$) of best arm identification in the
fixed-confidence (resp. fixed-budget) setting, as follows:
\[\kappa_{\text{C}}(\nu) = \inf_{\cA \ \delta-\text{PAC}} \limsup_{\delta
  \rightarrow 0} \frac{\bE_\nu[\tau]}{\log (1/\delta)}, \hspace{1cm}
\kappa_{\text{B}}(\nu)= \inf_{\cA \ \text{consistent}} \left(\limsup_{t
    \rightarrow \infty} -\frac{1}{t}\log p_t(\nu)\right)^{-1}.\]
Heuristically, for a given bandit model $\nu$ and a given $\delta>0$, a
fixed-confidence optimal strategy uses an average number of samples of order
$\kappa_{C}(\nu)\log (1/\delta)$, whereas a fixed-budget optimal strategy
uses approximately $t=\kappa_{B}(\nu)\log (1/\delta)$ draws in order to
ensure a probability of error at most equal to $\delta$. Most
of the existing performance bounds for the fixed confidence and fixed budget settings can be expressed using these complexity measures. 

The main goal of this paper is to determine $\kappa_{\text{C}}$ and
$\kappa_{\text{B}}$ for important classes of parametric bandit models, allowing
for a comparison between the fixed-confidence and fixed-budget settings.
Classical sequential testing theory provides a first element in that direction
in the simpler case of fully specified alternatives. Consider for instance the
case where $\nu_1$ and $\nu_2$ are Gaussian laws with the same known variance
$\sigma^2$, the means $\mu_1$ and $\mu_2$ being known up to a permutation. Denoting by $P$ the joint distribution of the paired samples $(X_s,Y_s)$, one
must choose between the hypotheses $H_0:P=\norm{\mu_1}{\sigma^2}\otimes
\norm{\mu_2}{\sigma^2}$ and $H_1:P=\norm{\mu_2}{\sigma^2}\otimes
\norm{\mu_1}{\sigma^2}$. It is known since \cite{Wald45SPRT} that among the
sequential tests such that type I and type II probabilities of error are both
smaller than $\delta$, the Sequential Probability Ratio Test (SPRT) minimizes
the expected number of required samples, and is such that
$\bE_\nu[\tau]=2\sigma^2/(\mu_1-\mu_2)^2\log(1/\delta)$. However, the batch test that minimizes both probabilities of error is the Likelihood Ratio test; it can be shown to require a sample size of order
$8\sigma^2/(\mu_1-\mu_2)^2\log(1/\delta)$ in order to ensure that both type I and type II error probabilities are smaller than $\delta$. Thus, when the sampling
strategy is uniform and the parameters are known, there is a clear gain in using
randomized stopping strategies. We show in the following that this conclusion
is not valid anymore when the values of $\mu_1$ and $\mu_2$ are not assumed to
be known.

\bigskip
\textbf{Related works.} Bandit models have received a considerable interest since their introduction by~\cite{Thompson33}.
An important focus was set on a different perspective, in which each observation is considered
as a reward: the agent aims at maximizing the cumulative rewards obtained until some horizon $t$. 
Equivalently, his goal is to  minimize the \emph{regret}
$R_t(\nu) = t\mu_{[1]} - \bE_\nu\left[\sum_{s=1}^{t}  Z_s\right]\;.$
Regret minimization, which is paradigmatic of the so-called \emph{exploration versus exploitation dilemma},
was introduced by \cite{Robbins:Freq52} and its complexity is well understood for parametric bandits. In generic one-parameter models,
\cite{LaiRobbins85bandits} prove that, with a proper notion of consistency adapted to regret minimization,  
\[\inf_{A \ \text{consistent}}\liminf_{t\rightarrow \infty}\frac{R_t(\nu)}{\log t} 
\geq \sum_{a : \mu_{a}<\mu_{[1]}} \frac{(\mu_{[1]} - \mu_a)}{\text{KL}(\nu_a,\nu_{[1]})}\;,\]
where $\text{KL}(\nu_i,\nu_j)$ denotes the Kullback-Leibler divergence between distributions $\nu_i$ and $\nu_j$.
Since then, non-asymptotic analyses of efficient algorithms matching this bound have been proposed. 
Optimal algorithms include the KL-UCB algorithm of \cite{KLUCB:Journal}---a variant of UCB1 (\cite{Auer:al02}) using 
informational upper bounds, Thompson Sampling (\cite{ALT12,AG:AISTAT13}), the DMED algorithm \citep{HondaTakemura11} and Bayes-UCB~\cite{AISTATS12}.
This paper is a first step in the attempt to similarly characterize the complexity of \emph{pure exploration}, where the goal is to determine the best arms without trying to maximize the cumulated observations.
 
The problem of best arm identification has received an important interest in
the 1950s as a particular case of 'ranking and identification problems'. The
literature on the subject goes beyond two-armed bandit models to finding the
$m>1$ best arms among $K>2$ arms, and sometimes introduces a relaxation
parameter $\epsilon>0$, such that arms within $\epsilon$ of the best arm should be
recommended. In the sequel, we always particularize the existing results to the
two-armed bandit models presented above. The first advances on this topic
are summarized in the monograph by \cite{Bechofer:al68}, who only consider the
fixed-confidence setting. In the same setting, algorithms introduced more
recently by
\cite{EvenDar:al06,Shivaram:al12,Gabillon:al12}
can be used to find the best arm in a two-armed bounded bandit model, in which
$\nu_1$ and $\nu_2$ are probability distributions on $[0,1]$.
Combining the upper bound on $\bE_\nu[\tau]$ for the LUCB algorithm of
\cite{Shivaram:al12} with the lower bound following from the work of
\cite{MannorTsi:04}, it can be shown that for bounded bandit models such that
$\mu_a\in[0, 1 - \alpha]$ for $a\in\{1,2\}$, there exists a constant $C_\alpha$
for which
\[{C_\alpha}/{(\mu_1-\mu_2)^{2}} \leq \kappa_C(\nu) \leq
{584}/(\mu_1-\mu_2)^{2}.\]

The fixed-budget setting has been studied recently by
\cite{Bubeck:BestArm10,Bubeck:alMult13}. In two-armed bandit problems, the
algorithms introduced in these papers boil down to sampling each arm $t/2$ times---$t$ denoting the total budget---and recommending the empirical best arm. A simple
upper bound on the probability of error of this strategy can be derived, and
this result paired with the lower bound of \cite{Bubeck:BestArm10} yields, for
bounded bandit models such that $\mu_a\in [\alpha;1-\alpha]$ for $a\in\{1,2\}$:
\[{(2/5)\alpha(1-\alpha)}/{(\mu_1-\mu_2)^{2}} \leq \kappa_B(\nu) \leq
{2}/(\mu_1-\mu_2)^{2}.\]

\cite{Bubeck:al11} show that in the fixed-budget setting any sampling
strategy designed to minimize regret performs poorly with respect to the
\textit{simple regret} $r_t := \mu^* - \mu_{\hat{S}_1}$, a quantity closely
related to the probability $p_t(\nu)$ of recommending the wrong arm. Therefore,
good strategies for best arm identification have to be quite different from
UCB-like strategies. We will show below that the complexities $\kappa_B(\nu)$ and $\kappa_C(\nu)$ of
pure-exploration involve information terms that are different from the Kullback-Leibler divergence featured in Lai and Robbins' lower bound on
regret.

\bigskip

\textbf{Contents of the paper.} Compared to existing results, we provide general lower bounds on
$\kappa_B(\nu)$ and $\kappa_C(\nu)$ that: \emph{(i)} are tighter, leading in
specific parametric cases to a precise evaluation of these complexities;
\emph{(ii)} do not require unnecessary support assumptions;
and {\emph(iii)} are stated in terms of information divergences between the
distributions $\nu_1$ and $\nu_2$ rather than in terms of the gap
$\mu_1-\mu_2$. As can be expected, we will indeed confirm that the inverse of
the squared gap $(\mu_1-\mu_2)^2$ is the relevant measure of complexity
only in the Gaussian case, and an approximation (in the spirit of Pinsker's inequality) for sub-Gaussian distributions.

Lower bounds on the sample complexity (resp. probability of error) of
algorithms using the uniform sampling strategy in the fixed-confidence
(resp. fixed-budget) setting are also derived and we show that for Gaussian
bandit models with different variances, there is a significant gain in using a
non-uniform sampling strategy. For Bernoulli bandits however, we show that
little can be gained by departing from uniform sampling, and we therefore
propose close-to-optimal tests both for the batch and sequential settings.  For
Gaussian bandits with a known common variance the optimal algorithm uses
uniform sampling. In this specific case, we propose an improved
$\delta$-PAC stopping rule, illustrating its performance through numerical
experiments. 

Our contributions follow from two main mathematical results: Lemma
\ref{lem:Cornerstone} provides a general relation between the expected
number of draws and Kullback-Leibler divergences of the arms' distributions, which is the key element to
derive the lower bounds. Lemma \ref{thm:subgaussian} is a tight deviation
inequality for martingales with sub-Gaussian increments, in the spirit of the
Law of Iterated Logarithm.
  
The paper is structured as follows. Section \ref{sec:GeneralBounds} presents a distribution-dependent lower bound on both $\kappa_B(\nu)$ and $\kappa_C(\nu)$ under the some identifiability assumption, as well as lower bounds for algorithms using
uniform sampling.  Gaussian bandit models are then studied in details in
Section \ref{sec:GaussianBandits}, and Bernoulli bandit models in Section
\ref{sec:Bernoulli}. Section \ref{sec:Experiments} includes a practical
illustration of the performance of matching algorithms for Gaussian bandits, as
well as a practical comparison of the fixed-confidence and fixed-budget
settings. The most important elements of proof are gathered in Section
\ref{sec:annex}, with the rest of the proofs in the Appendix.

\section{Lower Bounding the Complexity \label{sec:GeneralBounds}}

Introducing the Kullback-Leibler divergence of any two probability
distributions $p$ and $q$:
$$\K(p,q)  = 
\left\{
  \begin{array}{l}
    \int \log \left[\frac{dp}{dq}(x)\right]dp(x) \ \text{if} \ q \ll p,\\
    + \infty \ \text{otherwise},
  \end{array}
\right.
$$
we make the assumption that there exists a set $\cN$ such that for all
$\nu=(\nu_1,\nu_2)\in \cM$, for $a\in\{1,2\}, \nu_a\in\cN$ and that $\cN$
satisfies
\[
\forall p,q\in\cN, \ p\neq q \ \Rightarrow \ 0< \K(p,q) < + \infty.
\]
A class $\cM$ of bandit models satisfying this property is called
identifiable. For $\cM$ an identifiable class of bandit models, Theorem
\ref{thm:2arms} provides lower bounds on $\kappa_B(\nu)$ and $\kappa_C(\nu)$
for every $\nu\in\cM$. The proof of this theorem is based on changes of
distribution and detailed in Section \ref{sec:annex}.

\begin{theorem}\label{thm:2arms}
  Let $\mathbf{\nu}=(\nu_{1},\nu_{2})$ be a two-armed bandit model such that
  $\mu_1>\mu_2$.  In the fixed-budget setting, any consistent algorithm
  satisfies
  \[
  \limsup_{t\rightarrow \infty}-\frac{1}{t} \log p_t(\nu) \leq c^*(\nu),\ \ \
  \text{where} \ \ \ c^*(\nu) := \inf_{(\nu_1',\nu_2')\in\cM: \mu'_1 <\mu'_2}
  \max\left\{\K(\nu_1',\nu_{1}),\K(\nu_2',\nu_{2})\right\}.
  \]
  In the fixed-confidence setting any algorithm that is $\delta$-PAC on $\cM$ satisfies, when
  $\delta \leq 0.15$,
  \[
  \bE_\nu[\tau]\geq \frac{1}{c_*(\nu)} \log\left(\frac{1}{2 \delta}\right), \
  \ \ \text{where} \ \ \ c_*(\nu) := \inf_{(\nu_1',\nu_2')\in\cM:
    \mu_1'<\mu_2'}
  \max\left\{\K(\nu_{1},\nu_{1}'),\K(\nu_{2},\nu_{2}')\right\}.
  \]
\end{theorem}

In particular, Theorem \ref{thm:2arms} implies that $\kappa_B(\nu) \geq
1/{c^*(\nu)}$ and $\kappa_C(\nu) \geq
1/{c_*(\nu)}$. Proceeding similarly, one can obtain lower bounds for the algorithms
that use uniform sampling of both arms. The proof of the following result is easily
adapted from that of Theorem \ref{thm:2arms} (cf. Section~\ref{sec:annex}), using
that each arm is drawn $\tau/2$ times.
 
\begin{theorem}\label{thm:BoundUniform}Let $\mathbf{\nu}=(\nu_{1},\nu_{2})$ be a two-armed bandit model such that $\mu_1>\mu_2$. In the fixed-budget setting, any consistent algorithm 
  using a uniform sampling strategy satisfies
$$\limsup_{t\rightarrow\infty} -\frac{1}{t}\log p_t(\nu) \leq I^*(\nu) \ \ \ \text{where} \ \ \ I^*(\nu):=\inf_{(\nu_1',\nu_2')\in\cM: \mu'_1 <\mu'_2}\frac{\K\left(\nu_1',\nu_1\right)+\K\left(\nu_2',\nu_2\right)}{2}.$$
In the fixed-confidence setting, any algorithm that is $\delta$-PAC on $\cM$
and uses a uniform sampling strategy satisfies, for $\delta\leq 0.15$,
$$\hspace{1.5cm}\bE_\nu[\tau] \geq \frac{1}{I_*(\nu)}\log \frac{1}{2\delta} \ \ \ \text{where} \ \ \ I_*(\nu):=\inf_{(\nu_1',\nu_2')\in\cM: \mu'_1 <\mu'_2}\frac{\K\left(\nu_1,\nu_1'\right)+\K\left(\nu_2,\nu_2'\right)}{2}.$$
\end{theorem}

Obviously, one always has $I^*(\nu) \leq c^*(\nu)$ and $I_*(\nu) \leq c_*(\nu)$ suggesting
that uniform sampling can be sub-optimal. It is possible to give explicit expressions for the quantities
$c^*(\nu),c_*(\nu)$ and $I^*(\nu),I_*(\nu)$ for specific classes of parametric
bandit models that will be considered in the rest of the paper. In the
case of Gaussian bandits with known variance (see Section~\ref{sec:GaussianBandits}):
\begin{equation}\cM = \{ \nu =
  \left(\norm{\mu_1}{\sigma_1^2},\norm{\mu_2}{\sigma_2^2}\right) :
  (\mu_1,\mu_2)\in\R^2, \mu_1\neq \mu_2\},\label{set:Gaussian}\end{equation}
one obtains
\begin{eqnarray*}
  c^*(\nu)=c_*(\nu) = \frac{(\mu_1-\mu_2)^2}{2(\sigma_1+\sigma_2)^2} \ \ \ \text{and} \ \ \ I^*(\nu)=I_*(\nu)=\frac{(\mu_1 - \mu_2)^2}{4(\sigma_1^2 + \sigma_2^2)}.
\end{eqnarray*}
Hence, the lower bounds of Theorem \ref{thm:2arms} are equal in this case, and
we provide in Section \ref{sec:GaussianBandits} matching upper bounds confirming
that indeed $\kappa_B(\nu) = \kappa_C(\nu)$. In addition, the observation that
$2 I^*(\nu) \geq c^*(\nu) \geq I^*(\nu)$ implies that, except when $\sigma_1 =
\sigma_2$, strategies based on uniform sampling are sub-optimal.

The values of $c^*(\nu)$ and $c_*(\nu)$ can also be computed for 
\textit{canonical one-parameter exponential families} with density with respect to
some reference measure given by
\begin{equation}f_\theta(x)=A(x)\exp(\theta x - b(\theta)), \ \ \text{for} \
  \theta\in \Theta \subset \R.\label{DensityExpo}\end{equation}
We consider the class of bandit models
\[\cM = \{ \nu = (\nu_{\theta_1},\nu_{\theta_2}) :
(\theta_1,\theta_2)\in\Theta^2, \theta_1\neq \theta_2\}\] where
$\nu_{\theta_a}$ has density $f_{\theta_a}$ given by (\ref{DensityExpo}).
Using the shorthand $\Kb(\theta,\theta')=\K(\nu_{\theta},\nu_{\theta'})$ for
$(\theta,\theta')\in \Theta^2$, one can show that, for $\nu$ such that
$\mu_1>\mu_2$: 
\[
\begin{array}{ccccl}
  c^*(\nu) &= &  \inf_{\theta\in\Theta} \max \left(\Kb(\theta,\theta_1),\Kb(\theta,\theta_2)\right)
  & = & \Kb(\theta^*,\theta_1),  \ \ \text{where} \ \ \Kb(\theta^*,\theta_1)=\Kb(\theta^*,\theta_2),\\
  c_*(\nu) &=& \inf_{\theta\in\Theta} \max \left(\Kb(\theta_1,\theta),\Kb(\theta_2,\theta)\right)
  & = & \Kb(\theta_1,\theta_*),  \ \ \text{where} \ \ \Kb(\theta_1,\theta_*)=\Kb(\theta_2,\theta_*). 
\end{array}
\]
The coefficient $c^*(\nu)$ is known as the \textit{Chernoff information}
$\Kb^*(\theta_1,\theta_2)$ between the distributions $\nu_{\theta_1}$ and
$\nu_{\theta_2}$ (see \cite{Cover:Thomas} and \cite{COLT13} for earlier notice
of the relevance of this quantity in the best arm selection problem).  By
analogy, we will also denote $c_*(\nu)$ by
$\Kb_*(\theta_1,\theta_2)=\Kb(\theta_1,\theta_*)$.

For exponential family bandits the quantities $c^*(\nu)$ and $c_*(\nu)$ are not
equal in general, although it can be shown that it is the case when the
log-partition function $b(\theta)$ is (Fenchel) self-conjugate (e.g., for Gaussian and
exponential variables). In Section~\ref{sec:Bernoulli}, we will focus on the
case of Bernoulli models for which $c^*(\nu) > c_*(\nu)$. By exhibiting a
matching strategy in the fixed-budget case, we will show that this implies that
$\kappa_C(\nu)>\kappa_B(\nu)$ in this case.

\section{The Gaussian Case \label{sec:GaussianBandits}} 
We study in this Section the class of two-armed Gaussian bandit models with
known variances defined by (\ref{set:Gaussian}), where $\sigma_1$ and
$\sigma_2$ are fixed.
In this case, we observed above that the lower bounds of Theorem \ref{thm:2arms}
are similar, because $c^*(\nu) = c_*(\nu)$. We prove in this section that indeed
\[\kappa_C(\nu)=\kappa_B(\nu)=\frac{2(\sigma_1+\sigma_2)^2}{(\mu_1-\mu_2)^2}\]
by exhibiting strategies that reach these performance bounds. These strategies
are based on the simple recommendation of the empirical best arm but use
non-uniform sampling in cases where $\sigma_1$ and $\sigma_2$ differ. When
$\sigma_1 = \sigma_2$ we provide in Theorem~\ref{thm:PACFC} an improved
stopping rule that is $\delta$-PAC but results in a significant reduction
of the running time of fixed-confidence tests.

\subsection{Fixed-Budget Setting}

We consider the simple family of \textit{static strategies} that draw $n_1$
samples from arm 1 followed by $n_2=t-n_1$ samples of arm 2, and then choose
arm 1 if $\hat{\mu}_{1,n_1}<\hat{\mu}_{2,n_2}$, where $\hat{\mu}_{i,n_i}$ denotes the
empirical mean of the $n_i$ samples from arm $i$. Assume for instance that
$\mu_1>\mu_2$. Since $\hat{\mu}_{1,n_1} - \hat{\mu}_{2,n_2}-\mu_1 + \mu_2 \sim
\norm{0}{{\sigma_1^2}/{n_1}+{\sigma_2^2}/{n_2}}$, the probability of error of
such a strategy is easily upper bounded as:
\begin{eqnarray*}
  \bP\left(\hat{\mu}_{1,n_1} < \hat{\mu}_{2,n_2}\right) &
  \leq & 
  \exp\left(-\left(\frac{\sigma_1^2}{n_1}+\frac{\sigma_2^2}{n_2}\right)^{-1}\frac{(\mu_1-\mu_2)^2}{2}\right). 
\end{eqnarray*}
The right hand side is minimized when $n_1/(n_1+n_2) = {\sigma_1}/{(\sigma_1 +
  \sigma_2)}$, and the static strategy drawing $n_1 = \left\lceil\sigma_1
  t/(\sigma_1 + \sigma_2) \right\rceil$ times arm 1 is such that
\[ \liminf_{t \rightarrow \infty} -\frac{1}{t}\log p_t(\nu) \geq \frac{(\mu_1 -
  \mu_2)^2}{2(\sigma_1 + \sigma_2)^2}\;,\] which matches the bound of Theorem
\ref{thm:2arms} for Gaussian bandit models.

\subsection{Fixed-Confidence Setting}

\subsubsection{Equal Variances \label{subsec:equalvar}}

We start with the simpler case $\sigma_1=\sigma_2=\sigma$, where the quantity $I_*(\nu)$ introduced in Theorem
\ref{thm:BoundUniform} coincides with $c_*(\nu)$, which suggests that uniform sampling could be optimal.
A uniform sampling strategy is equivalent to collecting paired samples
$(X_s,Y_s)$ from both arms. The difference $X_s - Y_s$ is Gaussian with
mean $\mu=\mu_1-\mu_2$ and a $\delta$-PAC algorithm is equivalent to a sequential
test of $H_0:\mu<0$ versus $H_1:\mu>0$ such that the probability of error is
uniformly bounded by $\delta$. \cite{Robbins70LIL} proposes such a test that
stops after a number of samples
\begin{equation}\label{Robbins}
  \tau = \inf\left\{ t\in 2\N^*: \left|\sum_{s=1}^{t/2} (X_s - Y_s)\right| > \sqrt{2\sigma^2t\beta(t,\delta)}\right\} 
  \ \text{with} \ \beta(t,\delta)=\frac{t+1}{t}\log\left(\frac{t+1}{2\delta}\right)
\end{equation}
and recommends the empirical best arm. This procedure belongs to the class of
\textit{elimination strategies}, introduced by
\cite{Jennison:al84} who derive a lower bound on the sample complexity of any
$\delta$-PAC \textit{elimination} strategy---whereas our lower bound applies to
\textit{any} $\delta$-PAC algorithm---matched by Robbins' algorithm, that is, $\lim_{\delta \rightarrow
  0}{\bE_\nu[\tau]}/{\log\frac{1}{\delta}} = {8\sigma^2}/{(\mu_1-\mu_2)^2}.$
Therefore, Robbins' rule (\ref{Robbins}) yields an optimal strategy matching
our general lower bound of Theorem \ref{thm:2arms} in the particular case of
Gaussian distributions with common known variance.

Note that any elimination strategy that is $\delta$-PAC and uses a threshold function smaller than Robbins' also 
matches our asymptotic lower bound, while being strictly more efficient than Robbins' rule. 
For practical purpose, it is therefore interesting to exhibit smaller \textit{exploration rates} $\beta(t,\delta)$ 
leading to a $\delta$-PAC algorithm. The probability of error of such an algorithm is upper bounded, for example for $\mu_1<\mu_2$ by
\begin{equation}\label{PACExplain}
  \bP_\nu\left(\exists k\in\N : \sum_{s=1}^{k} \frac{X_s - Y_s - (\mu_1-\mu_2)}{\sqrt{2\sigma^2}}
    > \sqrt{2k\beta(2k,\delta)}\right) =  \bP\left(\exists k\in\N : S_k > \sqrt{2k\beta(2k,\delta)}\right)
\end{equation}
where $S_k$ is a sum of $k$ i.i.d. variables of distribution
$\norm{0}{1}$. \cite{Robbins70LIL} obtains a non-explicit confidence region of
risk at most $\delta$ by choosing
$\beta(2k,\delta)=\log\left({\log(k)}/{\delta}\right) + o(\log\log(k))$. The
dependency in $k$ is in some sense optimal, because the Law of Iterated
Logarithm (LIL) states that
$\limsup_{k\rightarrow\infty}{S_k}/\sqrt{2k\log\log(k)} =1$ almost
surely. Recently, \cite{Jamieson:al13LILUCB} proposed an explicit confidence region
inspired by the LIL. However, Lemma~1 of \citep{Jamieson:al13LILUCB} cannot be used to upper bound (\ref{PACExplain}) by
$\delta$ and we provide in Section \ref{sec:annex} a result derived independently (Lemma \ref{thm:subgaussian}) 
that achieves this goal and yields the following result.

\begin{theorem} \label{thm:PACFC} For $\delta$ small enough, the elimination strategy (\ref{Robbins}) is $\delta$-PAC with
  \begin{equation}\label{choiceBetaBetter}\beta(t,\delta)=\log \frac{1}{\delta}
    + \frac{3}{4} \log\log \frac{1}{\delta} +
    \frac{3}{2}\log(1+\log(t/2)).\end{equation}
\end{theorem}

\subsubsection{Mismatched Variances\label{subsec:General}}

In the case where $\sigma_1\neq \sigma_2$, we rely on an $\alpha$-elimination strategy, described in Algorithm~\ref{AlgoBox:Elimination}. For $a=1,2$, $\hat{\mu}_{a}(t)$ denotes
the empirical mean of the samples gathered from arm $a$ up to time $t$. 
The algorithm is based on a non-uniform sampling strategy governed by the parameter $\alpha\in(0,1)$
which ensures that, at the end of every round $t$, $N_1(t)=\lceil \alpha t\rceil$, $N_2(t)=t -
\lceil \alpha t\rceil$ and $\hat{\mu}_1(t)-\hat{\mu}_2(t) \sim \norm{\mu_1 - \mu_2}{\sigma_t^2(\alpha)}$. The sampling schedule used here is thus deterministic.

\begin{algorithm}[t]
  \caption{$\alpha$-Elimination\label{AlgoBox:Elimination}}
  \begin{algorithmic}[1]
    \REQUIRE Exploration function $\beta(t,\delta)$, parameter $\alpha$.
    \STATE \textit{Initialization}: $\hat{\mu}_1(0)=\hat{\mu}_2(0)=0$, $\sigma^2_{0}(\alpha)=1$. $t=0$.
    \WHILE {$|\hat{\mu}_1(t) - \hat{\mu}_2(t)| \leq \sqrt{2\sigma^2_{t}(\alpha)\beta(t,\delta)}$}
    \STATE $t=t+1$.
    \STATE If $\lceil\alpha t \rceil = \lceil \alpha(t-1) \rceil$,  $A_{t}=2$, else $A_t=1$. 
    \STATE Observe $Z_t\sim \nu_{A_t}$ and compute the empirical means $\hat{\mu}_1(t)$ and $\hat{\mu}_2(t)$.
    \STATE Compute $\sigma_t^2(\alpha)=\sigma_1^2/\lceil\alpha t\rceil + \sigma_2^2/(t-\lceil \alpha t \rceil).$ 
    \ENDWHILE
    \RETURN $a=\argmax{a=1,2} \ \hat{\mu}_a(t)$
  \end{algorithmic}
\end{algorithm}

Theorem \ref{thm:MatchingFC} shows that the $\sigma_1/(\sigma_1 +
\sigma_2)$-elimination algorithm, with a suitable exploration rate, is
$\delta$-PAC and matches the lower bound on $\bE_\nu[\tau]$, at least
asymptotically when $\delta \rightarrow 0$. Its proof can be found in Appendix \ref{proof:MatchingFC}.

\begin{theorem}\label{thm:MatchingFC} If $\alpha=\sigma_1/(\sigma_1 + \sigma_2)$, the $\alpha$-elimination strategy using the exploration rate  
  $\beta(t,\delta)=\log\frac{t}{\delta} + 2 \log\log(6t)$ is $\delta$-PAC on
  $\cM$ and satisfies, for every $\nu\in\cM$, for every $\epsilon>0$,
$$\bE_\nu[\tau] \leq (1+\epsilon)\frac{2(\sigma_1 + \sigma_2)^2}{(\mu_1-\mu_2)^2}\log\left(\frac{1}{\delta}\right) + \underset{\delta \rightarrow 0}{o_{\epsilon}}\left(\log\left(\frac{1}{\delta}\right)\right).$$
\end{theorem}

\begin{remark} When $\sigma_1=\sigma_2$, $1/2$-elimination reduces, up to
  rounding effects, to the elimination procedure described in
  Section~\ref{subsec:equalvar}, for which Theorem
  \ref{thm:PACFC} suggests an exploration rate
  of order $\log(\log(t)/\delta)$. As the feasibility of this exploration rate
   when $\sigma_1\neq\sigma_2$ is yet to be established, we focus on Gaussian
  bandits with equal variances in the numerical experiments of Section~\ref{sec:Experiments}. 
\end{remark}

\section{The Bernoulli Case \label{sec:Bernoulli}}

We consider in this section the class of Bernoulli bandit models defined by
\[\cM = \{ \nu=\left(\cB(\mu_1),\cB(\mu_2)\right) : (\mu_1,\mu_2)\in ]0;1[^2,
\mu_1 \neq \mu_2 \},\] where each arm can be equivalently parametrized by the natural parameter of the exponential
family, $\theta_a = \log(
{\mu_a}/{(1-\mu_a)})$. Following the notation of Section \ref{sec:GeneralBounds}, the
Kullback-Leibler divergence between two Bernoulli distributions can be
either expressed as a function of the means,
$\K(\cB(\mu_1),\cB(\mu_2))$, or of the natural parameters,
$\Kb(\theta_1,\theta_2)$.

In this Section, we prove that $\kappa_C(\nu)>\kappa_B(\nu)$ for Bernoulli
bandit models (Proposition \ref{prop:Bernoulli2}). To do so, we first introduce
a static strategy matching the lower bound of Theorem \ref{thm:2arms} in the
fixed-budget case (Proposition \ref{prop:Bernoulli1}). This strategy is
reminiscent of the algorithm exhibited for Gaussian bandits in Section
\ref{sec:GaussianBandits} and uses parameter-dependent non uniform
sampling. This strategy is not directly helpful in practice but we show that it
can be closely approximated by an algorithm using the uniform sampling
strategy. In the fixed-confidence setting we similarly conjecture that little can be
gained from using a non-uniform sampling strategy and propose an algorithm based
on a non-trivial stopping strategy that is believed to match the bound of
Theorem \ref{thm:BoundUniform}.

\begin{proposition}\label{prop:Bernoulli1} Let $\alpha(\theta_1,\theta_2)$ be
  defined by
  \[\alpha(\theta_1,\theta_2)=\frac{\theta^*-\theta_1}{\theta_2-\theta_1} \ \ \
  \text{where} \ \ \Kb(\theta^*,\theta_1)=\Kb(\theta^*,\theta_2).\] 
  For all $t$, the static strategy that allocates $\left\lceil\alpha(\theta_1,\theta_2) t\right\rceil$   samples to arm  1
  , and recommends the empirical best arm, satisfies 
  $p_t(\nu) \leq \exp(-K^*(\theta_1,\theta_2)t)$.
\end{proposition}

This result, proved in Appendix \ref{proof:ConcExp}, shows in particular that for every $\nu\in\cM$ there exists a consistent
static strategy such that
\[\liminf_{t\rightarrow \infty} - \frac{1}{t}\log p_t \geq
\Kb^*(\theta_1,\theta_2), \ \ \ \text{and hence that} \ \
\kappa_B(\nu)=\frac{1}{\Kb^*(\theta_1,\theta_2)}.\] However, as
$\alpha(\theta_1,\theta_2)$ depends in the Bernoulli case on the unknown means of the arms,
this optimal static strategy is not useful in practice. 
So far, it is not known whether there exists a \textit{universal} strategy such that $p_t(\nu)\leq \exp(-\Kb^*(\theta_1,\theta_2)t)$ for all bandit instance $\nu\in\cM$.

For Bernoulli bandit models it can be checked that for all
$\nu\in\cM$, 
$\Kb_*(\theta_1,\theta_2)<\Kb^*(\theta_1,\theta_2)$. This fact together with
Proposition \ref{prop:Bernoulli1} and Theorem \ref{thm:2arms} yields the
following inequality.

\begin{proposition}\label{prop:Bernoulli2}
  For all $\nu\in\cM$, $\kappa_C(\nu)>\kappa_B(\nu)$.
\end{proposition}

\begin{figure}[t]
  \centering
  \
  \includegraphics[height=5cm]{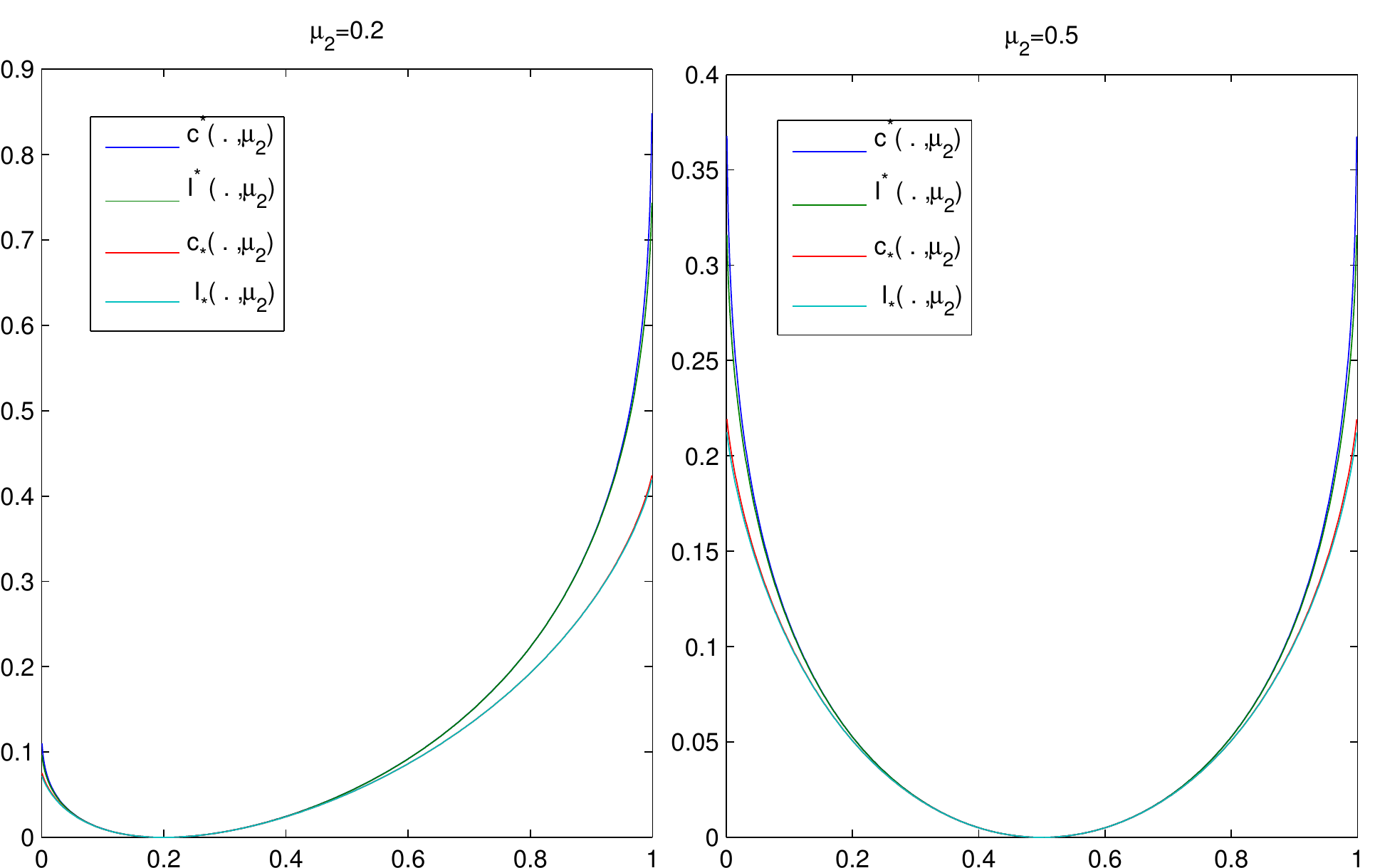}
  \caption{\label{fig:CompareComplexities} Comparison of different information terms in Bernoulli bandit models.}
\end{figure}

In the specific case of Bernoulli distributions, there is a strong incentive to use uniform sampling: the quantities $I^*(\nu)$ and
$I_*(\nu)$ introduced in Theorem \ref{thm:BoundUniform} appear to be very close to
$c^*(\nu)$ and $c_*(\nu)$ respectively. This fact is illustrated in Figure
\ref{fig:CompareComplexities}, on which we represent these different
quantities, that are functions of the means $\mu_1,\mu_2$ of the arms, as a
function of $\mu_1$, for two fixed values of $\mu_2$.  Therefore, algorithms
matching the bounds of Theorem \ref{thm:BoundUniform} provide upper bounds on
$\kappa_B(\nu)$ (resp. $\kappa_C(\nu)$) very close to $1/c^*(\nu)$
(resp. $1/c_*(\nu)$). In the fixed-budget setting, Lemma \ref{lem:ConcExp} shows that the strategy with uniform sampling that
recommends the empirical best arm, satisfies $p_t(\nu)\leq e^{-tI^*(\nu)}$, and
matches the bound of Theorem \ref{thm:BoundUniform} (see Remark \ref{rem:Match} in Appendix \ref{proof:ConcExp}). 
Hence, problem-dependent optimal strategy
described above can be approximated by a very simple, universal algorithm.
Similarly, finding an algorithm for the fixed-confidence setting sampling the
arms uniformly and matching the bound of Theorem \ref{thm:BoundUniform} is a
crucial matter. This boils down to finding a good stopping rule. In all the
algorithms studied so far, the stopping rule was based on the difference of the
empirical means of the arms. For Bernoulli arms, such a strategy can be
analyzed with the tools provided in this paper: the algorithm stopping for $t$
such that $\hat{\mu}_{1,t/2} - \hat{\mu}_{2,t/2} >
\sqrt{{2\beta(t,\delta)}/{t}}$ with $\beta(t,\delta)$ as in Theorem
\ref{thm:PACFC} is $\delta$-PAC and its expected running time bounded by
${2}/{(\mu_1-\mu_2)^2}\log\frac{1}{\delta} + o
\left(\log\frac{1}{\delta}\right)$. Yet, Pinsker's inequality implies that
$I_*(\mu_1,\mu_2)>(\mu_1-\mu_2)^2/2$ and this algorithm is thus not optimal with
respect to Theorem \ref{thm:BoundUniform}. The approximation $I_*(\mu_1,\mu_2)
= (\mu_1-\mu_2)^2/(8\mu_1(1-\mu_1)) + o\left((\mu_1-\mu_2)^2\right)$ suggests
that the loss with respect to the optimal error exponent is particularly
significant when both means are close to 0 or 1. The stopping rule we propose
to circumvent this drawback is
the following:
\begin{equation}\tau = \inf\left\{ t \in 2\N^* :
    tI_*(\hat{\mu}_{1,t/2},\hat{\mu}_{2,t/2}) > \log
    \left((\log(t)+1)/{\delta}\right)\right\}. \label{OptimalBernoulli}
\end{equation}
This algorithm is related to the KL-LUCB algorithm of \cite{COLT13}. Indeed,
$I_*(x,y)$ mostly coincides with $\Kb_*(\cB(x),\cB(y))$
(Figure~\ref{fig:CompareComplexities}) and a closer examination shows that the
stopping criterion in KL-LUCB for two arms is exactly of the form
$tK_*(\cB(\hat{\mu}_{1,t/2}),\cB(\hat{\mu}_{2,t/2})) > \beta(t,\delta)$. The
results of \cite{COLT13} show in particular that the algorithm based on
(\ref{OptimalBernoulli}) is provably $\delta$-PAC for appropriate choices
of $\beta(t,\delta)$. However, by analogy with the result of Theorem \ref{thm:PACFC}
we believe that the analysis of \cite{COLT13} is too conservative
and that the proposed approach should be $\delta$-PAC for exploration rates 
$\beta(t,\delta)$ that grow as a function of $t$ only at rate $\log\log t$.

\section{Numerical Experiments and Discussion \label{sec:Experiments}}

The goal of this Section is twofold: to compare results obtained in 
the fixed-budget and fixed-confidence settings and to illustrate the improvement
resulting from the adoption of the reduced exploration rate of Theorem \ref{thm:PACFC}.

\begin{figure}[t] \centering
  \vspace{-0.3cm}
  \begin{center}
    \mbox{ \hspace{-0.4cm}
      \includegraphics[height=5.9cm]{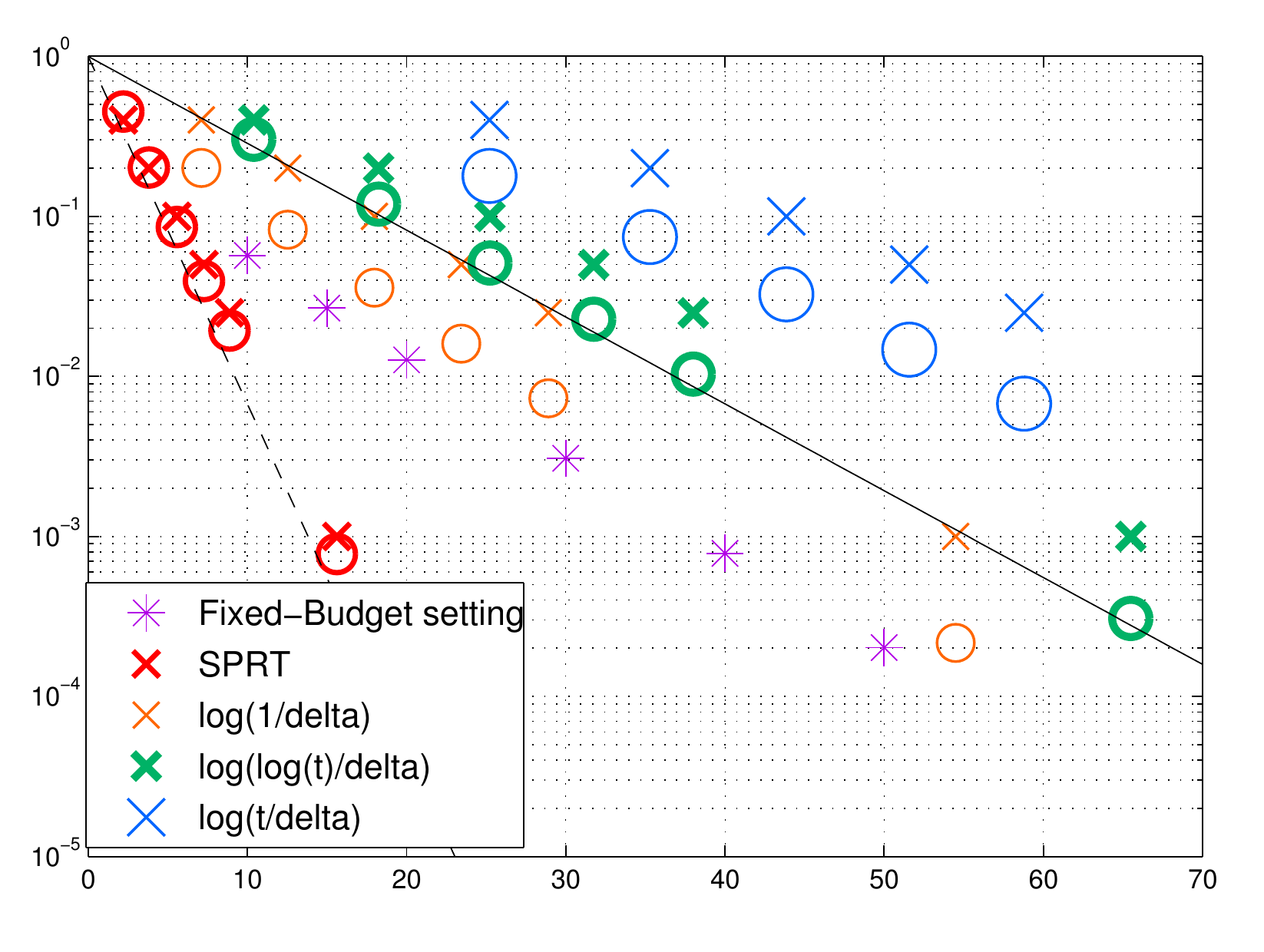}
      \hspace{-0.2cm}
      \includegraphics[height=5.9cm]{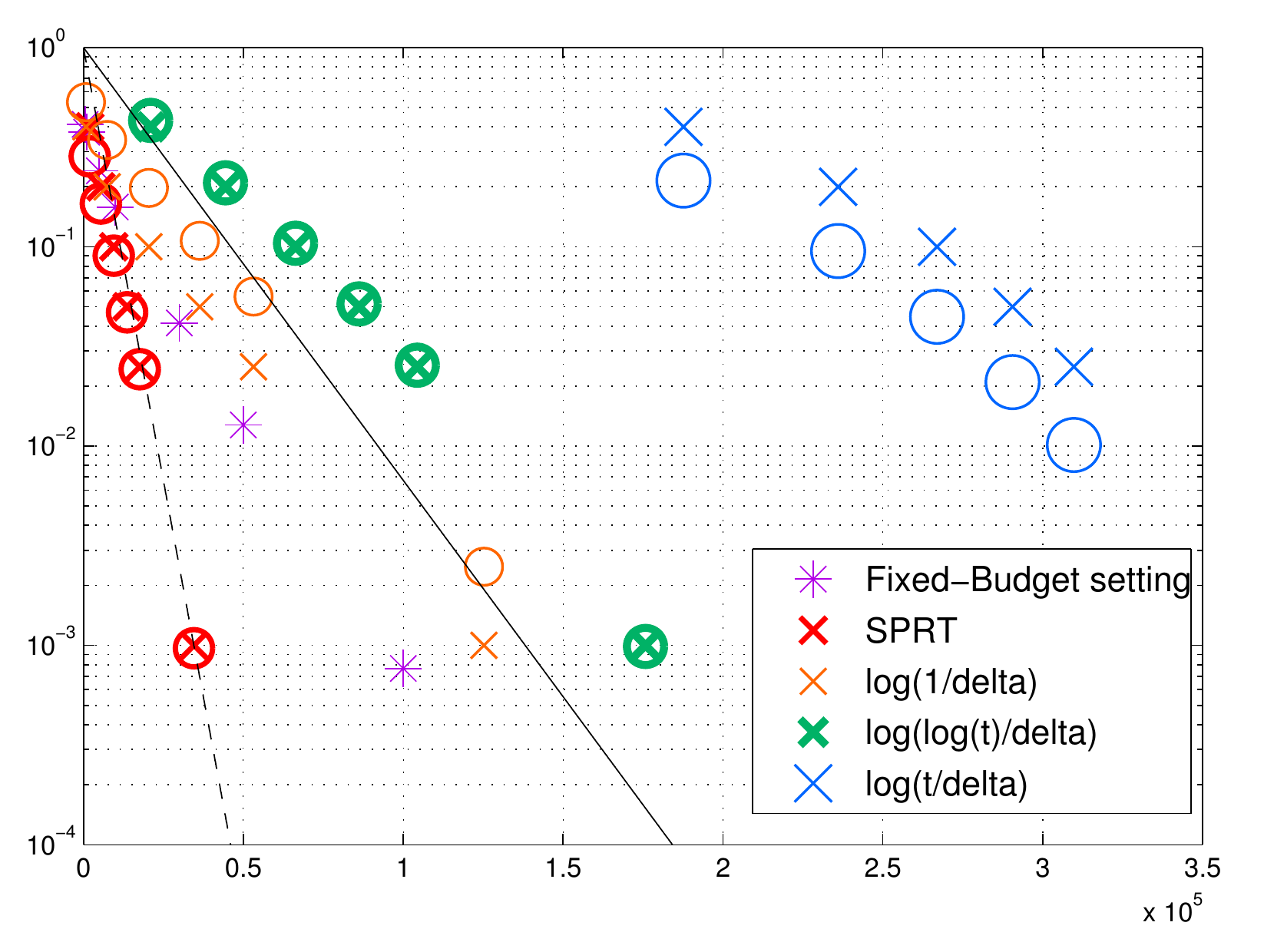}
}
\vspace{-0.4cm}
\caption{Experimental results (descriptions in text; plots best viewed in
  color).}
\label{fig:ResultsGaussian}
\vspace{-0.6cm}
\end{center}
\end{figure}

In Figure \ref{fig:ResultsGaussian}, we consider two bandit models: the 'easy'
one is $\norm{0.5}{0.25}\otimes \norm{0}{0.25}$, $\kappa=8$ (left) and the
'difficult' one is $\norm{0.01}{0.25}\otimes \norm{0}{0.25}$, $\kappa=20000$
(right). In the fixed-budget setting, stars ('*') report the probability of
error $p_n(\nu)$ as a function of $n$. In the
fixed-confidence setting, we plot both the empirical probability of error by circles
('O') and the specified maximal error probability $\delta$ by crosses ('X') as
a function of the empirical average of the running times. Note the logarithmic scale used
for the probabilities on the y-axis.  All results are
averaged on $N=10^6$ independent Monte Carlo replications. For comparison
purposes, a plain line represents the theoretical rate $x \mapsto \exp(-x/\kappa)$ which is a straight line on the log scale.

In the fixed-confidence setting, we report results for algorithms of the form
(\ref{Robbins}) with $g(t,\delta)=\sqrt{2\sigma^2t\beta(t,\delta)}$ for three
different exploration rates $\beta(t,\delta)$. The exploration rate we consider
are: the provably-PAC rate of Robbins' algorithm $\log({t}/{\delta})$ (large blue symbols),
the conjectured 'optimal' exploration rate
$\log({(\log(t)+1)}/{\delta})$, almost provably $\delta$-PAC according to
Theorem \ref{thm:PACFC} (bold green symbols), and the rate $\log ({1}/{\delta})$,
which would be appropriate if we were to perform the stopping test only at a
single pre-specified time (orange symbols).  For
each algorithm, the log probability of error is approximately a linear function
of the number of samples, with a slope close to $-1/\kappa$, where $\kappa$ is
the complexity. We can visualize the gain in sample complexity achieved by
smaller exploration rates, but while the rate $\log((\log(t)+1)/\delta)$
appears to guarantee the desired probability of error across all problems, the
use of $\log(1/\delta)$ seems too risky, as one can see that the probability of
error becomes larger than $\delta$ on difficult problems. To illustrate the
gain in sample complexity when the means of the arms are known, we add in red
the SPRT algorithm mentioned in the introduction along with the theoretical
relation between the probability of error and the expected number of samples,
materialized as a dashed line. The SPRT stops for $t$ such that
$|(\mu_1-\mu_2)(S_{1,t/2} - S_{2,t/2})|>\log(1/\delta)$.

Robbins' algorithm is $\delta$-PAC and matches the complexity (which is
illustrated by the slope of the measures), though in practice the use of the
exploration rate $\log((\log(t)+1)/\delta)$ leads to huge gain in terms of
number of samples used. It is important to keep in mind that running times
play the same role as error exponents and hence the threefold increase of average
running times observed on the rightmost plot of Figure~\ref{fig:ResultsGaussian} when
using $\beta(t,\delta)=\log({t}/{\delta})$ is really prohibitive.
This illustrates the asymptotic nature of our notion of
complexity: the leading term in $\bE_\nu[\tau]$ is indeed $\kappa
\log(1/\delta)$ but there is a second-order constant term which is not
negligible for fixed value of $\delta$. \cite{Jamieson:al13LILUCB} implicitly
consider an alternative complexity for Gaussian bandit models with common known
variance: for a fixed value of $\delta$, if $\Delta=\mu_1 -\mu_2$, they show
that when the gap $\Delta$ goes to zeros, the sample complexity if of order
some constant---that depends on $\delta$---multiplied by $\Delta^{-2}
\log\log \Delta^{-2}$.

If one compares on each problem the results for the fixed-budget setting to
those for the best $\delta$-PAC algorithm (in green), one can see that to
obtain the same probability of error, the fixed-confidence algorithm needs an
average number of samples of order at least twice larger than the deterministic
number of samples required by the fixed-budget setting algorithm. This remark
should be related to the fact that a $\delta$-PAC algorithm is designed to be
uniformly good across all problems, whereas consistency is a weak requirement
in the fixed-budget setting: any strategy that draws both arm infinitely often
and recommends the empirical best is
consistent. Figure~\ref{fig:ResultsGaussian} shows that when the values of
$\mu_1$ and $\mu_2$ are unknown, the sequential version of the test is no more
preferable to its batch counterpart and can even become much worse if the
exploration rate $\beta(t,\delta)$ is chosen too conservatively. This
observation should be mitigated by the fact that the sequential (or
fixed-confidence) approach is adaptive with respect to the difficulty of the
problem whereas it is impossible to predict the efficiency of a batch (or
fixed-budget) experiment without some prior knowledge regarding the problem under
consideration.

\section{Elements of Proof\label{sec:annex}}

\subsection{Proof of Theorem \ref{thm:2arms}}

The cornerstone of the proof of all the lower bounds given in this paper is Lemma~\ref{lem:Cornerstone}
which relates the probabilities of the same event under two different models to the expected
number of draws of each arm. Its proof, which may be found in Appendix \ref{sec:proofCorn}, 
encapsulates the technical aspects of the change of distributions. $N_a(t)$ denotes the number of draws
of arm $a$ up to round $t$ and $N_a=N_a(\tau)$ is the total number of draws of
arm $a$ by some algorithm $\cA$.

\begin{lemma} \label{lem:Cornerstone} Let $\nu$ and $\nu'$ be two bandit
  models.  For any $A\in\cF_\tau$ such that $0<\bP_\nu(A)<1$
  \begin{equation}
    \bE_\nu[N_1] \K(\nu_{1},\nu'_{1}) + \bE_\nu[N_2] \K(\nu_{2},\nu'_{2})\geq d(\bP_\nu(A),\bP_{\nu'}(A)),\label{ineq:Basic}
  \end{equation}
  where
  $d(x,y):=\K(\cB(x),\cB(y))=x\log (x/y) + (1-x)\log\big((1-x)/(1-y)\big).$
  
\end{lemma}

Without loss of generality, assume that the bandit model $\nu=(\nu_1,\nu_2)$ is
such that $a^*=1$. Consider any alternative bandit model $\nu'=(\nu_1',\nu_2')$
in which $a^*=2$ and the event $A=(\hat{a}=1)$ where $\hat{a}$ is the arm
chosen by algorithm $\cA$. Clearly $A\in\cF_\tau.$

\paragraph{Fixed-Confidence Setting.}

Let $\cA$ be a $\delta$-PAC algorithm. Then it is correct on both $\nu$ and
$\nu'$ and satisfies $\bP_\nu(A)\geq 1-\delta$ and $\bP_{\nu'}(A) \leq \delta$.
Using monotonicity properties of $d$ (for example $x\mapsto d(x,y)$ is increasing when $x>y$
and decreasing when $x<y$) and inequality (\ref{ineq:Basic}) in Lemma
\ref{lem:Cornerstone} yields $\bE_\nu[N_1]\K(\nu_{1},\nu_{1}') + \bE_\nu[N_2]\K(\nu_{2},\nu_{2}') \geq d(\delta,1-\delta)$, and hence
\begin{equation*}
  \bE_\nu[\tau] \geq \frac{d(\delta,1-\delta)}{\max_{a=1,2}\K(\nu_{a},\nu_{a}')},
\end{equation*}
using that $\tau = N_1 + N_2$. Optimizing over the possible model $\nu'$
satisfying $\mu_1'<\mu_2'$ to make the right hand side of the inequality as
large as possible gives the result, using moreover that for $\delta \leq 0.15$,
it can be shown that $d(1-\delta,\delta)\geq \log ({1}/(2\delta))$.

\paragraph{Fixed-Budget Setting.} 
Inequality
(\ref{ineq:Basic}) in Lemma \ref{lem:Cornerstone} applied to $A$ yields
\[ \bE_{\nu'}[N_1(t)] \K(\nu'_{1},\nu_{1}) + \bE_{\nu'}[N_2(t)] \K(\nu'_{2},\nu_{2})
\geq d(\bP_{\nu'}(A),\bP_{\nu}(A)).\] Note that $p_t(\nu)=1-\bP_{\nu}(A)$ and
$p_t(\nu')= \bP_{\nu'}(A)$.  As algorithm $\cA$ is correct on both $\nu$ and
$\nu'$, for every $\epsilon>0$ there exists $t_0(\epsilon)$ such that for all
$t\geq t_0(\epsilon)$, $\bP_{\nu'}(A) \leq \epsilon \leq \bP_{\nu}(A)$.  For
$t\geq t_0(\epsilon)$,
\[
\bE_{\nu'}[N_1(t)] \K(\nu'_{1},\nu_{1}) + \bE_{\nu'}[N_2(t)] \K(\nu'_{2},\nu_{2}) \geq
d(\epsilon,1-p_t(\nu)) \geq (1-\epsilon)\log \frac{1 - \epsilon}{p_t(\nu)} +
\epsilon \log {\epsilon}.
\]
Taking the limsup (denoted by $\overline{\lim}$) and letting $\epsilon$ go to zero, one can show that
\begin{eqnarray*}
\overline{\lim}-\frac{1}{t} \log p_t(\nu) &\leq& \overline{\lim} \Bigl[\frac{\bE_{\nu'}[N_1(t)]}{t}\K(\nu_{1}',\nu_{1}) + \frac{\bE_{\nu'}[N_2(t)]}{t}\K(\nu_{2}',\nu_{2})\Bigr]\leq \max_{a=1,2}\K(\nu_{a}',\nu_{a}). 
\end{eqnarray*}
Optimizing over the possible model $\nu'$ satisfying $\mu_1'<\mu_2'$ to make
the right hand side of the inequality as small as possible gives the result.

\subsection{Proof of Theorem \ref{thm:PACFC}}

According to (\ref{PACExplain}), the proof of Theorem \ref{thm:PACFC} boils
down to finding an exploration rate such that $\bP(\exists t\in \N^*:
S_t>\sqrt{2\sigma^2 t \beta(t,\delta)})\leq \delta$, where $S_t=X_1+\dots +
X_t$ is a sum of i.i.d. normal random variable. Lemma \ref{thm:subgaussian}
provides such a confidence region. Its proof can be found in Appendix
\ref{proof:DevIneq}.

\begin{lemma}\label{thm:subgaussian}
  Let $\zeta(u) = \sum_{k\geq 1} k^{-u}$.  Let $X_1, X_2,\dots$ be independent
  random variables such that, for all $\lambda\in\R$,
  $\phi(\lambda):=\log\bE[\exp(\lambda X_1)] \leq \lambda^2\sigma^2/2$.  For
  every positive integer $t$ let $S_t = X_1 + \dots + X_t$.  Then, for all
  $\beta>1$ and $x\geq \frac{8}{(e-1)^2}$,
  \[\bP\Bigl(\exists t\in \N^* : {S_t} > \sqrt{2\sigma^2t(x +
      \beta\log\log(et))}\Bigr) \leq \sqrt{e}\, \zeta\Bigl(\beta
    \bigl(1-\frac{1}{2x}\bigr)\Bigr)\Bigl(\frac{\sqrt{x}}{2\sqrt{2}}+1\Bigr)^{\beta}
  \exp(-x).\]
\end{lemma}

Let $\beta(t,\delta)$ be of the form $\beta(t,\delta)=\log \frac{1}{\delta} + c
\log\log \frac{1}{\delta} + d\log\log(et)$, for some constants $c>0$ and $d>1$.  Lemma
\ref{thm:subgaussian} yields
\[\bP\Bigl(\exists t \in \N : S_t > \sqrt{2\sigma^2 t
    \beta(t,\delta)}\Bigr)\leq \zeta\Bigl(d\bigl(1-\frac{1}{2(z+c\log
      z)}\bigr)\Bigr)\frac{\sqrt{e}}{(2\sqrt{2})^d} \frac{(\sqrt{z + c\log z}
  + \sqrt{8})^d}{z^c} \delta,\] where $z:=\log \frac{1}{\delta}>0$. To
upper bound the above probability by $\delta$, at least for large values of $z$ (which corresponds to small values of $\delta$), it suffices to choose the
parameters $c$ and $d$ such that 
\[\sqrt{e}\,\zeta\Bigl(d\bigl(1-\frac{1}{2(z+c\log z)}\bigr)\Bigr)
\frac{1}{(2\sqrt{2})^d}\frac{(\sqrt{z + c\log z} + 2\sqrt{2})^d}{z^c} \leq 1.\]
For $c={d}/{2}$, the left hand side tends to
${\sqrt{e}\zeta\left(d\right)}/{(2\sqrt{2})^d}$ when $z$ goes to infinity,
which is smaller than 1 for $d\geq 1.47$. Thus, for $\delta$ small enough, the desired inequality
holds for $d = {3}/{2}$ and $c={3}/{4}$, which corresponds to the
exploration rate of Theorem \ref{thm:PACFC}.

\section{Conclusion}

We provide distribution-dependent lower bounds for best-arm identification in
the context of two-armed bandit models. 
These bounds involve information-theoretic quantities that reflect the typical causes of failure, which are different from those appearing in regret analysis.
For Gaussian and Bernoulli bandit models, we exhibit matching algorithms showing that these bounds are (mostly) tight,
highlighting relationships between the complexities of the fixed-budget and
fixed-confidence settings. Our numerical experiments illustrate the
significance of using appropriate exploration rates in the context of best
arm(s) identification and we believe that Lemma \ref{lem:Cornerstone} can be
adapted to deal with more general $K$-armed bandit scenarios.

These results suggest three practical implications for A/B testing. First, for
Binary and Gaussian-like responses with matched variances it is reasonable to
consider only tests---i.e., strategies using uniform sampling---rather than
general sequential sampling strategies.
Second, using a sequential stopping rule in this context is mostly of interest
because it does not requires prior knowledge of the complexity of the
problem. It should however not be expected to reduce the (average) running time
of the experiment for a given probability of error.
This leads to the third message regarding the utmost importance of using proper
(i.e., provably $\delta$-PAC but not too conservative) exploration rates when using
a sequential stopping rule.

\acks{We thank S\'ebastien Bubeck for fruitful discussions during the visit of the first author at Princeton University. This work was supported by the ANR-2010-COSI-002 grant of the French National Research Agency.}

\bibliography{kaufmann14biblio}

\appendix

\section{Proof of Lemma \ref{lem:Cornerstone}: Changes of Distributions\label{sec:proofCorn}}

Under the identifiability assumption, there exists a common measure $\lambda$
such that for all $\nu=(\nu_1,\nu_2)$, for all $a\in\{1,2\}$ $\nu_a$ has a
density $f_a$ with respect to $\lambda$.

Let $\nu\in\cM$ be a bandit model, and consider an alternative bandit model
$\nu'\in \cM$. $f_a,f_a'$ are the densities of $\nu_a,\nu_a'$ respectively and
one can introduce the log-likelihood ratio of the observations up to time $t$
under an algorithm $\cA$:
\[L_t :=
\sum_{s=1}^{t}\ind_{(A_s=1)}\log\left(\frac{f_{1}(Z_s)}{f'_1(Z_s)}\right) +
\sum_{s=1}^{t}\ind_{(A_s=2)}\log\left(\frac{f_{2}(Z_s)}{f'_2(Z_s)}\right).\]
The key element in a change of distribution is the following classical lemma
(whose proof is omitted) that relates the probabilities of an event under
$\bP_\nu$ and $\bP_{\nu'}$ through the log-likelihood ratio of the
observations. Such a result has often been used in the bandit literature for
$\nu$ and $\nu'$ that differ just from one arm (either $\nu_1=\nu'_1$ or
$\nu_2=\nu'_2$), for which the expression of the log-likelihood ratio is
simpler. As we will see, here we consider more general changes of
distributions.

\begin{lemma}\label{Garivier} Let $\sigma$ be any stopping time with respect to $\mathcal{F}_t$.  For every event $A\in\mathcal{F}_\sigma$ 
  (i.e. $A$ such that $A\cap (\sigma=t) \in \mathcal{F}_t$),
$$\bP_{\nu'}(A)=\bE_\nu[\ind_{A}\exp(-L_{\sigma})]$$
\end{lemma}

Let $\tau$ be the stopping rule of algorithm $\cA$. We start by showing that for all $A\in \cF_\tau$, $\bP_\nu(A)=0$ if and only if
$\bP_{\nu'}(A)=0$. Thus, if $0<\bP_\nu(A)<1$ one also has $0<\bP_{\nu'}(A)<1$
and the quantity $d(\bP_\nu(A),\bP_{\nu'}(A))$ in the statement of Lemma
\ref{lem:Cornerstone} is well defined.
Let $A\in \cF_\tau$. Lemma \ref{Garivier} yields $\bP_{\nu'}(A) =
\bE_\nu[\ind_{A}\exp(-L_\tau)]$.  Thus $\bP_{\nu'}(A)=0$ implies 
$\ind_{A}\exp(-L_\tau)=0 \ \bP_\nu-a.s$.  As $\bP_\nu(\tau < + \infty)=1$,
$\bP_\nu(\exp(L_\tau) > 0)=1$ and $\bP_{\nu'}(A)=0 \Rightarrow \bP_\nu(A)=0$. A
similar reasoning yields $\bP_\nu(A)=0 \Rightarrow
\bP_{\nu'}(A)=0$. 

We now prove Lemma \ref{lem:Cornerstone}. Let $A\in\cF_\tau$ be such that
$0<\bP_\nu(A)<1$. Then $0<\bP_{\nu'}(A)<1$. Lemma \ref{Garivier} and the
conditional Jensen inequality leads to
\begin{align*}
  \bP_{\nu'}(A) & = \bE_\nu[\exp(-L_\tau)\ind_{A}] =  \bE_\nu[\bE_\nu[\exp(-L_\tau) | \ind_{A}] \ind_{A}] \\
  & \geq \bE_\nu[\exp\left(-\bE_\nu[L_\tau | \ind_{A}]\right) \ind_{A}]  =  \bE_\nu[\exp\left(-\bE_\nu[L_\tau |A]\right) \ind_{A}] \\
  & = \exp\left(-\bE_\nu[L_\tau | A]\right) \bP_\nu(A).
\end{align*}
Writing the same for the event $\overline{A}$ yields
$\bP_{\nu'}(\overline{A})\geq \exp\left(-\bE_\nu[L_\tau | \overline{A}]\right)
\bP_\nu(\overline{A})$ and finally
$$\bE_\nu[L_\tau | A] \geq \log \frac{\bP_\nu(A)}{\bP_{\nu'}(A)} \ \ \ \text{and} \ \ \ \bE_\nu[L_\tau | \overline{A}] \geq \log \frac{\bP_\nu(\overline{A})}{\bP_{\nu'}(\overline{A})}.$$
Therefore one can write
\begin{eqnarray}
  \bE_\nu[L_\tau ] & = & \bE_\nu[L_\tau | A]\bP_\nu(A) + \bE_\nu[L_\tau | \overline{A}]\bP_\nu(\overline{A}) \nonumber \\
  & \geq & \bP_\nu(A) \log \frac{\bP_\nu(A)}{\bP_{\nu'}(A)} + \bP_\nu(\overline{A})\log \frac{\bP_\nu(\overline{A})}{\bP_{\nu'}(\overline{A})} = d(\bP_\nu(A),\bP_{\nu'}(A)) \label{LBLike}.
\end{eqnarray}
Introducing $(Y_{a,t})$, the sequence of i.i.d. samples successively observed from arm $a$, the log-likelihood ratio $L_t$ can be rewritten 
\[L_t = \sum_{a=1}^{2}\sum_{t=1}^{N_a(t)}\log\left(\frac{f_a(Y_{a,t})}{f'_a(Y_{a,t})}\right); \ \ \ \text{and} \ \ \ 
\bE_\nu\left[\log\left(\frac{f_a(Y_{a,t})}{f'_a(Y_{a,t})}\right)\right]=\K(\nu_{a},\nu'_{a}).\]
Applying Wald's Lemma (see for example \cite{Siegmund:SeqAn}) to $L_\tau =
\sum_{a=1}^{2}\sum_{t=1}^{N_a}\log\left(\frac{f_a(Y_{a,t})}{f'_a(Y_{a,t})}\right)$,
where $N_a=N_a(\tau)$ is the total number of draws of arm $a$, yields
$$\bE_\nu[L_\tau] = \bE_\nu[N_1] \K(\nu_{1},\nu'_{1})+\bE_\nu[N_2] \K(\nu_{2},\nu'_{2}),$$
which concludes the proof together with inequality (\ref{LBLike}).


\section{Proof of Lemma \ref{thm:subgaussian}: An Optimal Confidence
  Region \label{proof:DevIneq}}

We start by stating three technical lemmas, whose proofs are partly omitted.

\begin{lemma}\label{lem:sandwich}
  For every $\eta>0$, every positive integer $k$, and every integer $t$ such
  that $(1+\eta)^{k-1} \leq t \leq (1+\eta)^k$,
  \[\sqrt{\frac{(1+\eta)^{k-1/2}}{t}} + \sqrt{\frac{t}{(1+\eta)^{k-1/2}}} \leq
  (1+\eta)^{1/4} + (1+\eta)^{-1/4}\;.\]
\end{lemma}

\begin{lemma}\label{lem:boundeta}
  For every $\eta>0$,
  \[A(\eta) := \frac{4}{\left((1+\eta)^{1/4} + (1+\eta)^{-1/4}\right)^2} \geq
  1-\frac{\eta^2}{16}.\]
\end{lemma}

\begin{lemma}\label{lem:approxslice}
  Let $t$ be such that $(1+\eta)^{k-1} \leq t \leq (1+\eta)^k$.  Then, if
  $\lambda = \sigma^{-1}\sqrt{2zA(\eta)/(1+\eta)^{k-1/2}}$,
  \[\sigma\sqrt{2z} \geq \frac{A(\eta)z}{\lambda \sqrt{t}} +
  \frac{\lambda\sigma^2\sqrt{t}}{2}\;.\]
\end{lemma}
\textbf{Proof:}
\[
\frac{A(\eta)z}{\lambda\sqrt{t}} + \frac{\lambda\sigma^2\sqrt{t}}{2} =\frac{\sigma\sqrt{2zA(\eta)}}{2}\left(\sqrt{\frac{(1+\eta)^{k-1/2}}{t}} + \sqrt{\frac{t}{(1+\eta)^{k-1/2}}}\right)\\
\leq \sigma\sqrt{2z}
\]
according to Lemma~\ref{lem:sandwich}.  \qed

An important fact is that for every $\lambda\in\R$, because the $X_i$ are
$\sigma$-subgaussian, $W_t = \exp(\lambda S_t-t\frac{\lambda^2\sigma^2}{2}))$
is a super-martingale, and thus, for every positive $u$,
\begin{equation}
  \bP\left(\bigcup_{t \geq 1} \left\{\lambda S_t-t\frac{\lambda^2\sigma^2}{2} > u \right\}\right) \leq \exp(-u). 
  \label{EquMartingale}
\end{equation}

Let $\eta\in]0, e-1]$ to be defined later, and let $T_k = \N\cap
\left[(1+\eta)^{k-1}, (1+\eta)^{k}\right[$.
\begin{align*}
  &\bP\left(\bigcup_{t\geq 1} \left\{\frac{S_t}{\sigma\sqrt{2t}} > \sqrt{x +
        \beta\log\log(et)} \right\}\right)
  \leq \sum_{k=1}^\infty \bP \left(\bigcup_{t\in T_k} \left\{\frac{S_t}{\sigma\sqrt{2t}} > \sqrt{x + \beta\log\log(et)} \right\}\right)\\
  & \hspace{4cm} \leq \sum_{k=1}^\infty \bP \left(\bigcup_{t\in T_k}
    \left\{\frac{S_t}{\sigma\sqrt{2t}} > \sqrt{x +
        \beta\log\left(k\log(1+\eta)\right)}\right\}\right)\;.
\end{align*}

We use that $\eta\leq e-1$ to obtain the last inequality since this condition
implies $$\log(\log(e(1+\eta)^{k-1})\geq\log(k\log(1+\eta)).$$

For a positive integer $k$, let $z_k = x + \beta\log\left(k\log(1+\eta)\right)$
and $\lambda_k = \sigma^{-1}\sqrt{2z_kA(\eta)/(1+\eta)^{k-1/2}}$.

Lemma~\ref{lem:approxslice} shows that for every $t\in T_k$,
\[\left\{\frac{S_t}{\sigma\sqrt{2t}} > \sqrt{z_k} \right\}
\subset \left\{\frac{S_t}{\sqrt{t}} > \frac{A(\eta)z_k}{\lambda_k \sqrt{t}} +
  \frac{\sigma^2\lambda_k\sqrt{t}}{2}\right\}\;.\]

Thus, by inequality (\ref{EquMartingale}),
\begin{align*}
  \bP \left(\bigcup_{t\in T_k} \left\{\frac{S_t}{\sigma\sqrt{2t}} >
      \sqrt{z_k}\right\}\right)
  &\leq \bP \left(\bigcup_{t\in T_k} \left\{\frac{S_t}{\sqrt{t}} > \frac{A(\eta)z_k}{\lambda_k \sqrt{t}} + \frac{\sigma^2\lambda_k\sqrt{t}}{2}\right\}\right)\\
  & = \bP \left(\bigcup_{t\in T_k} \left\{\lambda_k S_t - \frac{\sigma^2\lambda_k^2 t}{2}> A(\eta)z_k\right\}\right)\\
  &\leq \exp\left(-A(\eta) z_k\right) = \frac{\exp(-A(\eta) x)
  }{(k\log(1+\eta))^{\beta A(\eta)}}\;.
\end{align*}

One chooses $\eta^2 = 8/x$ for $x$ such that $x\geq \frac{8}{(e-1)^2}$ (which
ensures $\eta\leq e-1$). Using Lemma~\ref{lem:boundeta}, one obtains that
$\exp(-A(\eta) x) \leq \sqrt{e} \exp(-x)$. Moreover,
\[
\frac{1}{\log(1+\eta)} \leq \frac{1+\eta}{\eta} =
\frac{\sqrt{x}}{2\sqrt{2}}+1\;.
\]
Thus,
\[\bP \left(\bigcup_{t\in T_k} \left\{\frac{S_t}{\sigma\sqrt{2t}} >
    \sqrt{z_k}\right\}\right) \leq \frac{\sqrt{e}}{k^{\beta
    A(\eta)}}\left(\frac{\sqrt{x}}{2\sqrt{2}}+1\right)^{\beta A(\eta)} \! \! \exp(-x)
\leq \frac{\sqrt{e}}{k^{\beta
    A(\eta)}}\left(\frac{\sqrt{x}}{2\sqrt{2}}+1\right)^{\beta} \exp(-x)\] and
hence,
\begin{align*}
  \bP\left(\bigcup_{t\geq 1}\left\{\frac{S_t}{\sigma\sqrt{2t}} > \sqrt{x + \beta\log\log(et)}\right\}\right) &\leq \sqrt{e}\zeta\left(\beta A(\eta)\right)\left(\frac{\sqrt{x}}{2\sqrt{2}}+1\right)^{\beta A(\eta)}  \exp\left(-x\right)\\
  &\leq \sqrt{e}\zeta\left(\beta
    \left(1-\frac{1}{2x}\right)\right)\left(\frac{\sqrt{x}}{2\sqrt{2}}+1\right)^{\beta}
  \exp\left(-x\right) \;,\end{align*}
using the lower bound on $A(\eta)$ given in Lemma \ref{lem:boundeta} and the fact that $A(\eta)$ is upper bounded by 1.


\section{Proof of Theorem \ref{thm:MatchingFC} \label{proof:MatchingFC}}

Let $\alpha = {\sigma_1}/{(\sigma_1+\sigma_2)}$. We first prove that with the
exploration rate $\beta(t,\delta)=\log(t/\delta) + 2\log\log(6t)$ the algorithm
is $\delta$-PAC.  Assume that $\mu_1>\mu_2$ and recall $\tau=\inf\{t\in\N :
|d_t| > \sqrt{2\sigma^2_t(\alpha)\beta(t,\delta)}\}$.  The probability of error
of the $\alpha$-elimination strategy is upper bounded by
\begin{eqnarray*}
  \bP_\nu\left(d_\tau \leq -\sqrt{{2\sigma_\tau^2(\alpha) \beta(\tau,\delta)}} \right) &\leq& \bP_\nu\left(d_\tau-(\mu_1-\mu_2) \leq - \sqrt{{2\sigma_\tau^2(\alpha) \beta(\tau,\delta)}} \right) \\
  &\leq& \bP_\nu\left(\exists t\in \N^* : d_t - (\mu_1 - \mu_2) <  - \sqrt{{2\sigma_t^2(\alpha) \beta(t,\delta)}}\right)\\
  & \leq & \sum_{t=1}^{\infty} \exp\left(-\beta(t,\delta)\right),
\end{eqnarray*}
by an union bound and Chernoff bound applied to $d_t - (\mu_1 - \mu_2)\sim
\norm{0}{\sigma^2_t(\alpha)}$.  The choice of $\beta(t,\delta)$ mentioned
above ensures that the series in the right hand side is upper bounded by $\delta$, which
shows the algorithm is $\delta$-PAC:
\[\sum_{t=1}^{\infty} e^{-\beta(t,\delta)}\leq \delta \sum_{t=1}^{\infty} \frac{1}{t(\log(6t))^2}\leq {\delta}\left(\frac{1}{(\log 6)^2} + \int_{1}^{\infty}\frac{dt}{t(\log(6t))^2}\right)
={\delta}\left(\frac{1}{(\log 6)^2} + \frac{1}{\log(6)}\right)\leq \delta.
\]

To upper bound the expected sample complexity, we start by upper bounding the
probability that $\tau$ exceeds some deterministic time $T$:
\begin{eqnarray*}
  \bP_\nu(\tau \geq T) &\leq& \bP_\nu\left(\forall t =1\dots T, \ d_t \leq \sqrt{{2\sigma_t^2(\alpha)\beta(t,\delta)}}\right) \leq \bP_\nu\left(d_T \leq \sqrt{{2\sigma_T^2(\alpha)\beta(T,\delta)}}\right) \\
  &=& \bP_\nu\left( d_T - (\mu_1 - \mu_2) \leq -\left[(\mu_1 - \mu_2)-\sqrt{{2\sigma_T^2(\alpha)\beta(T,\delta)}}\right]\right) \\
  &\leq & \exp\left(-\frac{1}{2\sigma_T^2(\alpha)}\left[(\mu_1 - \mu_2)-\sqrt{{2\sigma_T^2(\alpha)\beta(T,\delta)}}\right]^2\right).
\end{eqnarray*}
The last inequality follows from Chernoff bound and holds for $T$ such that
\begin{equation}(\mu_1 -
  \mu_2)>\sqrt{2\sigma^2_T(\alpha)\beta(T,\delta)}.\label{cdt}\end{equation}
Now, for $\gamma \in ]0,1[$ we introduce
\begin{eqnarray*}
  T^*_{\gamma} & := & \inf\left\{t_0 \in \N : \forall t \geq t_0, (\mu_1-\mu_2) - \sqrt{2\sigma^2_t(\alpha)\beta(t,\delta)} > \gamma (\mu_1 - \mu_2)\right\}. 
\end{eqnarray*}
This quantity is well defined as $ \sigma^2_t(\alpha)\beta(t,\delta)
\underset{t\rightarrow \infty}{\rightarrow} 0$. We then upper bound the
expectation of $\tau$:
\begin{eqnarray*}
  \bE_\nu[\tau]& \leq & T_\gamma^* + \! \sum_{T=T_{\gamma}^*+1}\bP\left(\tau \geq T\right) \\
  &\leq & T_\gamma^* + \! \sum_{T=T_{\gamma}^*+1}\exp\left(-\frac{1}{2\sigma_T^2(\alpha)}\left[(\mu_1 - \mu_2)-\sqrt{{2\sigma_T^2(\alpha)\beta(T,\delta)}}\right]^2\right) \\
  & \leq & T_\gamma^* + \! \sum_{T=T_{\gamma}^*+1}^\infty\exp\left(-\frac{1}{2\sigma_T^2(\alpha)}\gamma^2(\mu_1 - \mu_2)^2\right). 
\end{eqnarray*}
For all $t\in\N^*$, it is easy to show that the following upper bound on
$\sigma_t^2(\alpha)$ holds:
\begin{equation}
  \forall t\in \N, \  \sigma_t^2(\alpha) \leq \frac{(\sigma_1 + \sigma_2)^2}{t}\times \frac{t - \frac{\sigma_1}{\sigma_2}}{t-\frac{\sigma_1}{\sigma_2} -1}.
  \label{ineq:BoundSigma}
\end{equation}
Using the bound (\ref{ineq:BoundSigma}), one has
\begin{eqnarray*}
  \bE_\nu[\tau]& \leq & T_\gamma^* + \int_{0}^{\infty} \exp\left(-\frac{t}{2(\sigma_1 + \sigma_2)^2} \frac{t - \frac{\sigma_1}{\sigma_2}-1}{t - \frac{\sigma_1}{\sigma_2}}\gamma^2(\mu_1-\mu_2)^2\right)dt\\
  & = &  T_\gamma^* +  \frac{2(\sigma_1+\sigma_2)^2}{\gamma^2(\mu_1-\mu_2)^2}\exp\left(\frac{\gamma^2(\mu_1-\mu_2)^2}{2(\sigma_1+\sigma_2)^2}\right).
\end{eqnarray*}
We now give an upper bound on $T_\gamma^*$. Let $r\in[0,e/2-1]$. There exists $N_0(r)$
such that for $t\geq N_0(r)$, $\beta(t,\delta)\leq \log ({t^{1+r}}/{\delta})$.
Using also (\ref{ineq:BoundSigma}), one gets
$T_\gamma^*=\max(N_0(t),\tilde{T}_\gamma)$, where
$$\tilde{T}_\gamma = \inf \left\{ t_0 \in \N : \forall t \geq t_0, \frac{(\mu_1 - \mu_2)^2}{2(\sigma_1+\sigma_2)^2}(1-\gamma)^2 t > \frac{t - \frac{\sigma_1}{\sigma_2}-1}{t - \frac{\sigma_1}{\sigma_2}} \log \frac{t^{1+r}}{\delta}\right\}.$$
If $t > (1 + \gamma \frac{\sigma_1}{\sigma_2})/{\gamma}$ one has $(t -
\frac{\sigma_1}{\sigma_2}-1)/(t - \frac{\sigma_1}{\sigma_2}) \leq
{(1-\gamma)^{-1}}$. Thus $\tilde{T}_\gamma =\max((1 + \gamma
\frac{\sigma_1}{\sigma_2})/{\gamma}, T_\gamma')$, with
$$T_\gamma' = \inf \left\{ t_0 \in \N : \forall t\geq t_0, \exp\left(\frac{(\mu_1 - \mu_2)^2}{2(\sigma_1+\sigma_2)^2}(1-\gamma)^3 t\right) \geq \frac{t^{1+r}}{\delta}\right\}.$$ 
The following Lemma, whose proof can be found below, helps us bound this last
quantity.
\begin{lemma}\label{lem:Garivier2} For every $\beta,\eta>0$ and  $s\in[0,e/2]$, the following implication is true:
$$x_0 = \frac{s}{\beta}\log\left(\frac{e\log\left({1}/{(\beta^s\eta)}\right)}{\beta^s\eta}\right)  \ \ \ \Rightarrow \ \ \ \forall x\geq x_0, \ \ e^{\beta x} \geq \frac{x^s}{\eta}.$$
\end{lemma}
Applying Lemma \ref{lem:Garivier2} with $\eta=\delta$, $s=1+r$ and
\[\beta =\frac{(1-{\gamma})^3(\mu_1-\mu_2)^2}{2(\sigma_1+\sigma_2)^2}\]
leads to
\[T_\gamma' \leq\frac{(1+r)}{(1-\gamma)^3}\times \frac{2(\sigma_1 + \sigma_2)^2}{(\mu_1 - \mu_2)^2}\left[\log \frac{1}{\delta} + \log\log \frac{1}{\delta} \right] + R(\mu_1,\mu_2,\sigma_1,\sigma_2,\gamma,r),\]
with 
\[R(\mu_1,\mu_2,\sigma_1,\sigma_2,\gamma,r)=\frac{1+r}{(1-\gamma)^3}\frac{2(\sigma_1+\sigma_2)^2}{(\mu_1-\mu_2)^2}\left[1 + (1+r)\log\left(\frac{2(\sigma_1+\sigma_2)^2}{(1-\gamma)^3(\mu_1-\mu_2)^2}\right)\right].\]
Now for $\epsilon>0$ fixed, choosing $r$ and $\gamma$ small enough leads to
$$\bE_\nu[\tau] \leq (1+\epsilon)\frac{2(\sigma_1 + \sigma_2)^2}{(\mu_1 - \mu_2)^2}\left[\log \frac{1}{\delta} + \log\log\frac{1}{\delta}\right] + \cC(\mu_1,\mu_2,\sigma_1,\sigma_2,\epsilon),$$
where $\cC$ is a constant independent of $\delta$. This concludes the proof.

\qed

\paragraph{Proof of Lemma \ref{lem:Garivier2}} Lemma \ref{lem:Garivier2} easily
follows from the fact that for any $s,\eta>0$,
\[x_0 = s
  \log\left(\frac{e\log\left(\frac{1}{\eta}\right)}{\eta}\right) \ \
  \Rightarrow \ \ \ \forall x\geq x_0, \ \  e^x \geq \frac{x^s}{\eta}\]
Indeed, it suffices to apply this statement to $x=x\beta$ and $\eta=\eta\beta^s$. The mapping $x\mapsto e^x - x^s/\eta$ is increasing when $x\geq s$. As $x_0\geq s$,
it suffices to prove that $x_0$ defined above satisfies $e^{x_0}\geq x_0^s/\eta$. 
\begin{eqnarray*}
  \log\left(\frac{x_0^s}{\eta}\right) & = & s\log\left(s\log\left(\frac{e\log\frac{1}{\eta}}{\eta}\right)\right) + \log\frac{1}{\eta}  \\
  & = & s\left(\log(s) + \log \left[\log\frac{1}{\eta} + \log\left(e\log\frac{1}{\eta}\right)\right]\right) +  \log\frac{1}{\eta} \\
  &\leq&  s\left(\log(s) + \log \left[2\log\frac{1}{\eta}\right]\right)  +  \log\frac{1}{\eta}
\end{eqnarray*}
where we use that for all $y$, $\log(y) \leq \frac{1}{e}y$. Then
\begin{eqnarray*}
  \log\left(\frac{x_0^s}{\eta}\right) & \leq & s\left(\log(s) + \log(2) + \log\log\frac{1}{\eta} +  \log\frac{1}{\eta}\right). 
\end{eqnarray*}
For $s\leq \frac{e}{2}$, $\log(s)+\log(2) \leq 1$, hence
\begin{eqnarray*}
  \log\left(\frac{x_0^s}{\eta}\right) & \leq & s\left(1+ \log\log\frac{1}{\eta} +  \log\frac{1}{\eta}\right) = s \log\left(\frac{e\log\left(\frac{1}{\eta}\right)}{\eta}\right) = x_0,
\end{eqnarray*} 
which is equivalent to $e^{x_0} \geq \frac{x_0^s}{\eta}$ and concludes the proof.


\section{An Optimal Static Strategy for Bernoulli Bandit
  Models \label{proof:ConcExp}}

Bounding the probability of error of a static strategy using $n_1$ samples from
arm 1 and $n_2$ samples from arm 2 relies on the following Lemma, that applies
more generally to exponential families.

\begin{lemma} \label{lem:ConcExp} Let $(X_{1,t})_{t\in\N}$ and
  $(X_{2,t})_{t\in\N}$ be two independent i.i.d sequences, such that $X_{1,1}
  \sim \nu_{\theta_1}$ and $X_{2,1}\sim \nu_{\theta_2}$ belong to an
  exponential family. Assume that $\mu(\theta_1) > \mu(\theta_2)$. Then
  \begin{equation}\bP\left(\frac{1}{n_1}\sum_{t=1}^{n_1}X_{1,t}<
      \frac{1}{n_2}\sum_{t=1}^{n_2}X_{2,t}\right) \leq \exp(-(n_1 + n_2)
    g_\alpha(\theta_1,\theta_2)),\label{MainIneq}\end{equation}
  where $\alpha = \frac{n_1}{n_1 + n_2}$ and $g_{\alpha}(\theta_1,\theta_2)  := \alpha \Kb(\alpha\theta_1 + (1-\alpha)\theta_2,\theta_1) + (1-\alpha)\Kb(\alpha \theta_1 + (1-\alpha)\theta_2 , \theta_2).$
\end{lemma}
The function $\alpha \mapsto g_\alpha(\theta_1,\theta_2)$, can be maximized
analytically, and the value $\alpha^*$ that realizes the maximum is given by
\begin{eqnarray*}
  \Kb(\alpha^*\theta_1 + (1-\alpha^*)\theta_2,\theta_1) & = & \Kb(\alpha^*\theta_1 + (1-\alpha^*)\theta_2,\theta_2) \\
  \alpha^* \theta_1 + (1- \alpha^*)\theta_1 & = & \theta^* \\
  \alpha^* &=&\frac{\theta^* - \theta_2}{\theta_1 - \theta_2}
\end{eqnarray*}
where $\theta^*$ is defined by
$\Kb(\theta^*,\theta_1)=\Kb(\theta^*,\theta_2)=\Kb^*(\theta_1,\theta_2)$.  More
interestingly, the associated rate is such that
\[
g_{\alpha^*}(\theta_1,\theta_2)=\alpha^*\Kb(\theta^*,\theta_1) +
(1-\alpha^*)\Kb(\theta^*,\theta_2) =\Kb^*(\theta_1,\theta_2),
\]
which leads to Proposition \ref{prop:Bernoulli1}.

\begin{remark}\label{rem:Match} When $\mu_1>\mu_2$, applying Lemma \ref{lem:ConcExp} with $n_1=n_2=t/2$ yields 
\[\bP\left(\hat{\mu}_{1,t/2}<\mu_{2,t/2}\right) \leq \exp\left(-\frac{\Kb\left(\theta_1,\frac{\theta_1 +\theta_2}{2}\right)
+\Kb\left(\theta_2,\frac{\theta_1 +\theta_2}{2}\right)}{2}\,t\right) = \exp\big(-I_*(\nu)t\big),\]
which shows that uniform sampling matches the lower bound of Theorem \ref{thm:BoundUniform}.
\end{remark}

\paragraph{Proof of Lemma \ref{lem:ConcExp}}

The i.i.d. sequences $(X_{1,t})_{t\in\N}$ and $(X_{2,t})_{t\in\N}$ have respective
densities $f_{\theta_1}$ and $f_{\theta_2}$ where $f_\theta(x)=A(x)\exp(\theta
x - b(\theta)$ and $\mu(\theta_1)=\mu_1, \mu(\theta_2)=\mu_2$. $\alpha$ is such
that $n_1=\alpha n $ and $n_2=(1-\alpha)n$. One can write
\[
\bP\left(\frac{1}{n_1}\sum_{t=1}^{n_1}X_{1,t}-
  \frac{1}{n_2}\sum_{t=1}^{n_2}X_{2,t}< 0\right) = \bP\left(\alpha
  \sum_{t=1}^{n_2}X_{2,t} - (1-\alpha)\sum_{t=1}^{n_1}X_{1,t} \geq 0 \right).
\]
For every $\lambda > 0$, multiplying by $\lambda$, taking the exponential of the two sides and using
Markov's inequality (this technique is often referred to as Chernoff's method), one gets
\begin{multline*}
  \bP\left(\frac{1}{n_1}\sum_{t=1}^{n_1}X_{1,t}-
    \frac{1}{n_2}\sum_{t=1}^{n_2}X_{2,t}< 0\right) \leq \left(\bE_\nu[e^{\lambda \alpha
      X_{2,1}}]\right)^{(1-\alpha) n}
  \left(\bE_\nu[e^{\lambda (1-\alpha )X_{1,1}}]\right)^{\alpha n} \\
   = \exp\biggl(\underbrace{\left[(1-\alpha) \phi_{X_{2,1}}(\lambda \alpha) + \alpha
        \phi_{X_{1,1}}(-(1-\alpha)\lambda)\right]}_{G_\alpha(\lambda)}n\biggr)
\end{multline*}
with $\phi_X(\lambda)=\log \bE_\nu[e^{\lambda X}]$ for any random variable
$X$. If $X \sim f_\theta$ a direct computation gives $\phi_X(\lambda)=b(\lambda
+ \theta) - b(\theta)$. Therefore the function $G_\alpha(\lambda)$ introduced above
 rewrites
\[G_\alpha(\lambda) = (1-\alpha)(b(\lambda \alpha + \theta_2) - b(\theta_2)) + \alpha(b(\theta_1 - (1-\alpha)\lambda) - b(\theta_1)).\]
Using that $b'(x)=\mu(x)$, we can compute the derivative of $G$ and see that this
function as a unique minimum in $\lambda^*$ given by
\begin{eqnarray*}
\mu(\theta_1-(1-\alpha)\lambda^*) & = & \mu(\theta_2 + \alpha\lambda^*) \\
\theta_1-(1-\alpha)\lambda^* & = & \theta_2 + \alpha\lambda^* \\
\lambda^* & = & \theta_1-\theta_2,
  \end{eqnarray*}
using that $\theta \mapsto \mu(\theta)$ is one-to-one. One can also show that
\begin{eqnarray*}
  G(\lambda^*)& = &(1-\alpha)[ b(\alpha \theta_1 +(1-\alpha)\theta_2) - b(\theta_2)] + \alpha[b(\alpha \theta_1 + (1-\alpha)\theta_2) - b(\theta_1)]. 
\end{eqnarray*}
Using the expression of the KL-divergence between $\nu_{\theta_1}$
and $\nu_{\theta_2}$ as a function of the natural parameters:
$\Kb(\theta_1,\theta_2)= \mu(\theta_1)(\theta_1 - \theta_2) -
  b(\theta_1) + b(\theta_2)$, one can also show that 
\begin{align*}
  & \alpha \Kb(\alpha \theta_1 + (1-\alpha)\theta_2 , \theta_1) \\
  & \hspace{2cm} =  -\alpha(1-\alpha)\mu(\alpha\theta_1 + (1-\alpha)\theta_2)(\theta_1 - \theta_2) + \alpha[-b(\alpha \theta_1 + (1-\alpha)\theta_2) + b(\theta_1)] \\
  & (1-\alpha) \Kb(\alpha \theta_1 + (1-\alpha)\theta_2 , \theta_2) \\
  & \hspace{2cm}= \alpha(1-\alpha)\mu(\alpha\theta_1 +
  (1-\alpha)\theta_2)(\theta_1 - \theta_2) + (1-\alpha)[-b(\alpha \theta_1 +
  (1-\alpha)\theta_2) + b(\theta_2)]
\end{align*}
Summing these two equalities leads to
$$G(\lambda^*) = -\left[\alpha \Kb(\alpha \theta_1 + (1-\alpha)\theta_2 , \theta_1)  +  (1-\alpha) \Kb(\alpha \theta_1 + (1-\alpha)\theta_2 , \theta_2)\right] = - g_\alpha(\theta_1,\theta_2). $$
Hence the inequality $\bP\left(\frac{1}{n_1}\sum_{t=1}^{n_1}X_{1,t}<
  \frac{1}{n_2}\sum_{t=1}^{n_2}X_{2,t}\right)\leq \exp(G(\lambda^*)n)$ is
exactly (\ref{MainIneq}).

\end{document}